\documentclass[12pt,a4paper]{article}

\usepackage{latexsym}
\usepackage{amssymb}
\usepackage{amsthm}
\usepackage{amsmath}

\theoremstyle{definition}

\newcommand{\D}{\ensuremath{{\cal D}}}
\renewcommand{\S}{\ensuremath{{\cal S}}}

\newcommand{\mb}[1]{\ensuremath{\mathbb{#1}}}
\newcommand{\N}{\mb{N}}

\newcommand{\R}{\mb{R}}
\newcommand{\C}{\mb{C}}

\newcommand{\cl}[1]{\ensuremath{\mathrm{cl}[#1]}}

\newcommand{\G}{\ensuremath{{\cal G}}}

\newcommand{\WF}{\mathrm{WF}}
\newcommand{\singsupp}{\mathrm{sing supp}}

\renewcommand{\d}{\ensuremath{\partial}}

\newfont{\bl}{msbm10 scaled \magstep2}

\newtheorem{thm}{Theorem}
\newtheorem{lemma}[thm]{Lemma}
\newtheorem{prop}[thm]{Proposition}
\newtheorem{defn}[thm]{Definition}

\newtheorem{rem}[thm]{Remark}
\newtheorem{ex}[thm]{Example}

\newcommand{\bem}[1]{\vadjust{\rlap{\kern\hsize\thinspace\vbox%
                       to0pt{\hbox{${}_\clubsuit${\small\tt #1}}\vss}}}}

\newcommand{\map}{\ensuremath{\rightarrow}}

\newcommand{\F}{\ensuremath{{\cal F}}}

\newcommand{\notmid}{\mid\kern-0.5em\not\kern0.5em}

\newcommand{\al}{\alpha}
\newcommand{\be}{\beta}

\newcommand{\eps}{\varepsilon}

\newcommand{\vphi}{\varphi}

\newcommand{\om}{\omega}

\newcommand{\supp}{\mathop{\mathrm{supp}}}

\newcommand{\ovl}[1]{\overline{#1}}

\newcommand{\Ga}{\Gamma}
\newcommand{\Sig}{\Sigma}

\begin{document}

\begin{center}
{\bf \Large Microlocal properties of basic operations in Colombeau algebras}
\vskip2mm
{\large  G. H\"ormann \footnote{Present address: Center for Wave Phenomena,
Colorado School of Mines, Golden, CO 80401-1887, USA}, M. Kunzinger}

{\it Institut f\"ur Mathematik, Universit\"at Wien, Strudlhofgasse 4, A-1090 Wien,
Austria, E-mail: Guenther.Hoermann@univie.ac.at, Michael.Kunzinger@univie.ac.at}
\end{center}

{\small {\bf Abstract} The Colombeau algebra of generalized functions allows 
to unrestrictedly carry out 
products of distributions. We analyze this operation from a microlocal point of 
view, deriving a general inclusion relation for wave front sets of products in 
the algebra. Furthermore, we give explicit examples showing that the given result 
is optimal, i.e.\ its assumptions cannot be weakened. Finally, we discuss the 
interrelation of these results with the concept of pullback under smooth maps.

}

\section{Introduction}
Algebras of generalized functions in the sense of J.\ F.\ Colombeau provide an efficient
tool for the treatment of nonlinear problems involving singularities
(cf., e.g. \cite{c1}, \cite{c2}, \cite{clec}, \cite{MO}, \cite{GHKO} and the literature
cited therein). In particular, unrestricted multiplication (as well as a host of more 
general nonlinear operations) of distributions can be carried out in Colombeau algebras.
Moreover, starting with \cite{MO}, regularity theory has been 
introduced into the Colombeau framework and was extended to microlocal analysis with
applications to propagation of singularities in 
\cite{DPS}, \cite{NPS98}, \cite{GH99} and \cite{HdH00}. 

In the present paper we study microlocal 
properties of multiplication of generalized functions as well as of related operations
(like pullback) in this setting. Since unlike in the case of intrinsic multiplication of 
distributions the formation of products in the Colombeau algebra is not subject to
regularity conditions (``favorable position of the wave front sets''), new effects
can (and will) occur. Apart from deriving general results on inclusion 
relations for wave front sets of products the emphasis of our presentation
will be on providing examples illustrating these new effects. At the same time,
the examples will demonstrate that the mentioned inclusion relations are optimal
in the sense that the assumptions made to derive them cannot be weakened.

Concerning notation and terminology we basically follow \cite{MO}. Thus by ${\cal A}_0(\R)$
we denote the space of test functions on $\R$ with unit integral. For $1\le q$, ${\cal A}_q(\R)$
is the subspace of ${\cal A}_0(\R)$ consisting of those elements whose moments up to order
$q$ vanish. For $n\ge 1$, ${\cal A}_q(\R^n)$ is the space of $n$-fold tensor products
$\phi^{(n)}:=\phi\otimes\dots\otimes\phi$ with $\phi\in {\cal A}_0(\R)$. If $\phi \in \D(\R^n)$ 
is any test function we set $\phi_\eps(x) = \phi(x/\eps)/\eps^n$. Then the basic building
blocks of the Colombeau algebra of generalized functions are defined as follows:

${\cal E}_M(\R^n)$ is the set of all maps $R: {\cal A}_0(\R^n)\times \R^n \to \C$ which are
smooth in $x$ and satisfy: $\forall K \subset\subset \R^n$ $\forall \alpha \in \N_0^n$
$\exists N\in \N$ $\forall \phi \in {\cal A}_N(\R^n)$ $\exists c>0$ $\exists \eta>0$:
\begin{equation}\label{EM}
\sup_{x\in K} |\partial^\alpha R(\phi_\eps,x)| \le c \eps^{-N} \qquad (0<\eps<\eta)\,.
\end{equation}
${\cal N}(\R^n)$ is the subset of ${\cal E}_M(\R^n)$ consisting of those $R$ satisfying: 
$\forall K \subset\subset \R^n$ $\forall \alpha \in \N_0^n$ $\forall q\in \N$
$\exists p\in \N$ $\forall \phi \in {\cal A}_p(\R^n)$  $\exists c>0$ $\exists \eta>0$:
\begin{equation}\label{N}
\sup_{x\in K} |\partial^\alpha R(\phi_\eps,x)| \le c \eps^{q} \qquad (0<\eps<\eta)\,.
\end{equation}
Then the Colombeau algebra $\G(\R^n)$ is defined as the quotient ${\cal E}_M(\R^n)/{\cal N}(\R^n)$.
We note that to characterize ${\cal N}$ as a subspace of ${\cal E}_M$ it would suffice to suppose 
(\ref{N}) only for $\alpha=0$ (see \cite{MG00}, Th.\ 13.1). For the definition of $\G(\Omega)$ for 
$\Omega\subseteq \R^n$ open we refer to \cite{MO}. The subalgebra of compactly supported elements of
$\G$ will be denoted by $\G_c$. 

The respective definitions for the space $\G_\tau = {\cal E}_{M,\tau}/{\cal N}_\tau$ of tempered
Colombeau functions read:

${\cal E}_{M,\tau}(\R^n)$ is the set of all maps $R: {\cal A}_0(\R^n)\times \R^n \to \C$ which are
smooth in $x$ and satisfy: $\forall \alpha \in \N_0^n$
$\exists N\in \N$ $\forall \phi \in {\cal A}_N(\R^n)$ $\exists c>0$ $\exists \eta>0$:
\begin{equation}\label{EMtau}
|\partial^\alpha R(\phi_\eps,x)| \le c(1+|x|)^N \eps^{-N} \qquad (x\in \R^n,\, 0<\eps<\eta)
\end{equation}
${\cal N}_\tau(\R^n)$ is the subset of ${\cal E}_{M,\tau}(\R^n)$ consisting of those $R$ satisfying: 
$\forall \alpha \in \N_0^n$ $\forall q\in \N$
$\exists p\in \N$  $\forall \phi \in {\cal A}_p(\R^n)$ $\exists c>0$ $\exists \eta>0$:
\begin{equation}\label{Ntau}
|\partial^\alpha R(\phi_\eps,x)| \le c (1+|x|)^N \eps^{q} \qquad (x\in \R^n,\, 0<\eps<\eta)
\end{equation}

The canonical embedding of $\D'$ into $\G$ resp.\ of $\S'$ into $\G_\tau$ will consistently be denoted
by $\iota$. 
Also, equivalence classes of elements $R$ of ${\cal E}_M$ resp.\ ${\cal E}_{M,\tau}$
will be written as $\cl{(R(\phi,\,.\,))_\phi}$.

\section{Basic definitions, a first example}

The starting point for regularity theory and microlocal analysis in Colombeau algebras of generalized 
functions was the introduction of the subalgebra $\G^\infty$ of $\G$ by M.\ Oberguggenberger in 
\cite{MO}. $\G^\infty$ consists of those elements of $\G$ displaying uniform $\eps$-growth
in all derivatives. By \cite{MO}, Th.\ 25.2, $\G^\infty\cap \D' = {\cal C}^\infty$, an identity
on which all further regularity theory is based. The analogous notion $\G_\tau^\infty$ for tempered
Colombeau functions was introduced in \cite{GH99} where it was also shown that $\G_\tau^\infty\cap
\S' = {\cal O}_M$ (\cite{GH99}, Th.\ 16). Here ${\cal O}_M$ denotes the space of smooth functions
with at most polynomial growth in each derivative. 

$U\in \G_c$ is an element of $\G^\infty$ iff its Fourier transform (with respect to any damping measure)
is rapidly decreasing (\cite{GH99}, Th.\ 18). Based on this observation the concept of wave front set
in $\G$, first introduced in \cite{DPS}, has been (equivalently) stated in \cite{GH99} along the
lines of \cite{H}, Sect.\ 8.1. Thus for $U\in \G_c$ by $\Sigma_g(U)$ we denote the cone (in 
$\R^n\setminus 0$) which is the complement of those points possessing open conic neighborhoods
on which the Fourier transform of $U$ is rapidly decreasing. This notion is again independent
of the damping measure used in the definition of Fourier transform in $\G_\tau$. Then for $U\in \G(\Omega)$
and $x_0\in \Omega$, the cone of irregular directions at $x_0$ is 
\begin{equation}\label{cid}
\Sigma_{g,x_0}(U) = \bigcap_{\vphi \in \D(\Omega),\, \vphi(x_0)\not= 0} \Sigma_g(\vphi U)  
\end{equation}
The wave front set of $U$ is given by
\begin{equation}\label{wfdef}
\WF(U) = \{(x,\xi)\in \Omega\times \R^n\setminus 0 \mid \xi\in \Sigma_{g,x}(U)\}\,.  
\end{equation}
Finally, we shall make use of the concept of characteristic set of a linear differential
operator (introduced in \cite{DPS} for the case of the special Colombeau algebra):
\begin{defn} {\it
Let $P = \sum_{|\alpha|\le m} a_\alpha(x) \d^\al$ be a linear differential operator on $\Omega$ 
with coefficients in $\G(\Omega)$. $(x_0,\xi_0)\in \Omega\times \R^n\setminus \{0\}$ is not in the
characteristic set of $P$ if there exists a neighborhood $V_{x_0}$ of $x_0$, a
neighborhood $\Gamma_{\xi_0}$ of $\xi_0$, some $r\in \R$ and some $m\in \N_0$ such that 
$\forall \phi\in {\cal A}_m(\R^n)$  $\exists \eta>0$ $\exists C>0$ with
\begin{equation}
  \label{char0}
  |P_m(\phi_\eps,x,\xi)| \ge C\eps^r |\xi|^m, \qquad x\in V_{x_0}, \ \xi\in \Gamma_{\xi_0},
\ \eps \in (0,\eta)\,.
\end{equation}
}
\end{defn}
The following example, which was first introduced in \cite{GM}, gives a first application
of these concepts and introduces some methods that will repeatedly be used in the following
sections.
\begin{ex} \label{firstex}
We want to calculate the wave front set of the solution to
  \begin{equation}
    \label{sys}
    \begin{array}{rcl}
      (\d_t + a \d_x)U &=& 0 \\
       U\mid_{t=0} = U_0
    \end{array}
  \end{equation}
where $a$ denotes a bounded generalized constant and $U_0$ is allowed to be singular. 
Thus let $(a_\eps)_{\eps>0}$ be such that $|a_\eps|$ is bounded and denote by $a$ the 
generalized constant with representative $a(\phi) = a_{d(\phi)}$ (where $d(\phi)$
is the diameter of the support of $\phi$). Then for any $U_0\in \G(\R)$
the solution $U$ of (\ref{sys}) is given by the class of $U(\phi^{(2)},x,t) = 
u_{0}(\phi,x - a(\phi) t)$. Denote by
$B$ the set of limit points for $\eps\to 0$ of $(a_\eps)_{\eps>0}$. 
As was noted in \cite{GM},
\begin{equation} \label{singsupp}
\singsupp(U) \subseteq S := \{(x,t)\mid \exists b\in B \mbox{ s.t. } x-bt\in K \}
\end{equation}
where $K = \singsupp(U_0)$. In fact, let $(x_0,t_0)\not \in S$. Then by differentiating 
$u_{0}(\phi_\eps,x-a(\phi_\eps) t)$ and employing the boundedness of $a$ it follows that it 
suffices to show $x-a_{\eps} t \not\in K$ for $(x,t)$ in a neighborhood of $(x_0,t_0)$ 
and $\eps$ small which is obviously satisfied.

To begin with, let us determine the wave front set of $U$ for the particular initial 
value $U_0 = \delta = \cl{(\phi)_{\phi\in {\cal A}_0(\R)}}$, so $K=\{0\}$. 
Let $(x_0,t_0) \in S$ and let $\psi\in \D(\R^2)$, $\psi(x_0,t_0)\not=0$.
Then
\begin{eqnarray*}
(\psi u(\phi^{(2)}_\eps,.))^\wedge(\xi,\tau) = \int e^{-i(\xi x + \tau t)} \psi(x,t)
\eps^{-1}\phi\left(\frac{x-a_{d(\phi)\eps} t}{\eps}\right)\,d(x,t) \\
= \int e^{-i (\xi(\eps x + a_{d(\phi)\eps} t) +\tau t)}\psi(\eps x +a_{d(\phi)\eps} t,t) 
\phi(x)\,dx dt
\end{eqnarray*}
Setting $(\xi,\tau) = \om (\xi_0,\tau_0)$ this equals
\begin{equation}
\label{ft}
\int e^{i\om f_\eps(x,t)} \psi(\eps x +a_{d(\phi)\eps} t,t) \phi(x) \,dxdt\,.
\end{equation}
where
$f_\eps(x,t) = -(a_{d(\phi)\eps} \xi_0  + \tau_0) t - \eps \xi_0 x$. By \cite{H}, 
Th.\ 7.7.1, for any $k\in \N$ we obtain a constant $C$ independent of $\eps$ such that
\begin{equation} \label{hoer}
\begin{array}{l}
|(\psi u(\phi^{(2)}_\eps,.)^\wedge(\xi,\tau)| \le C \om^{-k} \\ 
\cdot \sum_{|\al|\le k} \sup (|D^\al(
\psi(\eps x+a_{d(\phi)\eps} t,t) \phi(x))| |Df_\eps(x,t)|^{|\al|-2k})\,.
\end{array}
\end{equation}
Here $|Df_\eps(x,t)|^2 = \eps^2 |\xi_0|^2 + |\tau_0 + a_{d(\phi)\eps} \xi_0|^2$, which
remains bounded away from $0$ (uniformly in $\eps$) for $\tau_0 \not\in -B\xi_0$.
Then from (\ref{hoer}) we conclude that any such pair $(\xi_0,\tau_0)$ is not contained
in $\Sigma_{g,(x_0,t_0)}(U)$. Denoting by $\Gamma_B$ the cone $\{(\xi,\tau)\mid \exists
b\in B \mbox{ with } \tau = -b\xi\}$ we have shown that
\begin{equation}
  \label{wfincl1}
  \WF_g(U) \subseteq S\times \Gamma_B \,.
\end{equation}
Conversely, let $\tau_0 = -b\xi_0$ for some $b\in B$, fix $\phi\in {\cal A}_0(\R)$ 
and choose a sequence $\eps_k\to 0$
with $a_{d(\phi)\eps_k}\to b$. Then from (\ref{ft}) we have 
$(\psi u(\phi^{(2)}_\eps,.)^\wedge(\xi,\tau) \to \int \psi(bt,t)\,dt$ 
which for an appropriate choice of $\psi$ is nonzero. Thus
$(\psi u(\phi^{(2)}_\eps,.)^\wedge$ is not rapidly decreasing in the direction 
$(\xi_0,\tau_0)$, so we also obtain the reverse inclusion of (\ref{wfincl1}). Summing up,
\begin{equation}
  \label{wfequ}
  \WF_g(U) = S\times \Gamma_B\,.
\end{equation}
We note that inclusion (\ref{wfincl1}) can even be obtained for general initial data
in $\G(\R)$ by propagation of singularities: in fact, by Th.\ 4 of \cite{DPS} we have
\begin{equation}
  \label{prop}
  \WF_g U \subseteq \mbox{Char} P \cup \WF_g P(U) \quad (U \in \G(\Omega))
\end{equation}
for any linear differential operator with coefficients in $\G^\infty(\Omega)$. Since
the right hand side of (\ref{sys}) is $0$ it therefore remains to determine the 
characteristic directions of the operator $P = \d_t + a \d_x$. 
Fix $\phi \in {\cal A}_0(\R)$ and set  $P_{\eps}(x,t;\xi,\tau) = P(\phi^{(2)}_\eps,x,t;\xi,\tau)
= i(\tau + a_\eps \xi)$ (to simplify notation we assume $d(\phi) = 1$). 
To show that $(x_0,t_0,\xi_0,\tau_0)$ is non-characteristic it suffices to prove the existence
of a neighborhood $V$ of $(x_0,t_0)$, a conic neighborhood of $(\xi_0,\tau_0)$
and constants $r\in \R$ (independent of $\phi$), $C>0$  and $\eta>0$ such that
\begin{equation}
  \label{char}
  |P_\eps(x,t;\xi,\tau)| \ge C \eps^r (|\xi| + |\tau|) \quad ((x,t)\in V, (\xi,\tau)\in \Gamma,
  0 < \eps < \eta)
\end{equation}
Let $\tau_0\not \in -B\xi_0$ and suppose $\xi_0 \not=0$ to begin with. Then there exists 
some $c>0$ such that $|\tau_0/\xi_0 + a_\eps| \ge c$ for $\eps$ small. Let 
$$
\Gamma = \{(\xi,\tau) \mid |\frac{\tau}{\xi} - \frac{\tau_0}{\xi_0}| < \frac{c}{2}\}
$$
For $(\xi,\tau)\in \Gamma$ we get
\begin{eqnarray*}
|\tau + a_\eps\xi| \ge |\xi|(|\frac{\tau_0}{\xi_0} + a_\eps| - |\frac{\tau}{\xi} - 
\frac{\tau_0}{\xi_0}|) \ge \frac{c}{2} |\xi| \ge \tilde c (|\xi| + |\tau|)
\end{eqnarray*}
so $(\xi,\tau)$ is noncharacteristic. On the other hand, if $\xi_0=0$ then we
choose $c$ such that $c|a_\eps| < 1/2$ for all $\eps$ and set $\Gamma = \{(\xi,\tau)
\mid |\xi/\tau| < c\}$. Then again
\begin{eqnarray*}
|\tau + a_\eps\xi| \ge |\tau| (1 - |a_\eps|\frac{|\xi|}{|\tau|}) \ge \frac{1}{2}|\tau|
\ge \tilde c (|\xi| + |\tau|)
\end{eqnarray*}
for $(\xi,\tau) \in \Gamma$. Thus
$\mbox{Char} P \subseteq \R^2 \times \Gamma_B$ which by (\ref{singsupp}) and 
(\ref{prop}) implies (\ref{wfincl1}).
\end{ex}

\begin{rem}
In \cite{GM} it is shown that inclusion (\ref{singsupp}) --- and consequently also (\ref{wfincl1})
--- may be strict.
\end{rem}

\section{The wave front set of a product}

In order to give a concise presentation of the following results we first collect
a few facts on cones in $\R^n$ resp.\ $\R^n\setminus 0$.

\begin{lemma}{\it \label{joshi}   
\
\begin{enumerate}
\item[(i)] If $\Sig_1$, $\Sig_2$ are closed cones in $\R^n$
such that $\Sig_1 \cap \Sig_2 = \{0\}$ then $\exists \al > 0$:
\[
  |\xi-\eta| \geq \al |\eta| \qquad \forall\xi\in\Sig_1, \forall\eta\in\Sig_2
  \, .
\]
\end{enumerate}
Let $\Ga_1$, $\Ga_2$ be closed cones in $\R^n\setminus 0$ such that 
$0 \not\in \Ga_1+\Ga_2$. Then

\begin{enumerate}
\item[(ii)] 
$\ovl{\Ga_1 + \Ga_2}^{\R^n\setminus 0} = (\Ga_1 + \Ga_2) \cup \Ga_1 \cup \Ga_2$.
\item[(iii)] For any open conic neighborhood 
$W$ of $\Ga_1+\Ga_2$ in $\R^n\setminus 0$ one can choose 
open conic neighborhoods $W_1$, $W_2$ in $\R^n\setminus 0$ of $\Ga_1$, $\Ga_2$,  respectively 
such that $W_1 + W_2 \subseteq W$.
\end{enumerate}
}
\end{lemma}
\begin{rem}\leavevmode
\begin{enumerate}
\item the second assertion is false (in gen.) if $0 \in \Ga_1+\Ga_2$ by 
  the following example (in $\R^3$) due to M. Grosser: 
Let $K_1 = \{\lambda \cdot(-1,t,t^2) \mid 0\le t \le 1, \lambda \ge 0\}$,
$K_2 = \{\lambda \cdot(1,t,t^2) \mid 0\le t \le 1, \lambda \ge 0\}$. Then
$\Ga_i = K_i\setminus 0$ ($i=1,2$) are closed cones in $\R^3\setminus 0$.
The sequence $\Ga_1 + \Ga_2 \ni \xi_n = n\cdot (-1,1/n,1/n^2) + n\cdot
(1,1/n,1/n^2) = (0,2,2/n)$ tends to $(0,2,0)$ as $n\to \infty$. But
$(0,2,0) \not\in  (\Ga_1+\Ga_2) \cup \Ga_1 \cup \Ga_2$: first, $(0,2,0)
= \lambda(-1,t,t^2) + \mu(1,s,s^2)$ implies $\lambda=\mu$, $\lambda\not=0$
and $s=t=0$, yielding $0=2$ in the second component, so $(0,2,0)\not\in 
\Ga_1+\Ga_2$. Also, $(0,2,0) \not\in \Ga_i$ $(i=1,2)$ 
by construction.

\item the third assertion is false (in gen.) if $0 \in \Ga_1+\Ga_2$ by 
the following example in $\R^2$: Let $\Ga_1 = \{(x,0) \mid x > 0\}$,
$\Ga_2 = \{(x,0) \mid x < 0\}$. Then the sum of any two open conic 
neighborhoods of $\Ga_1$, $\Ga_2$ is $\R^2$.
\end{enumerate}
\end{rem}
\begin{proof} \ \\
(i):
Otherwise there would be sequences $\xi_j\in\Sig_1$ and
$\eta_j\in\Sig_2$ ($j\in\N$) such that $|\xi_j-\eta_j| < |\eta_j|/j$ for all
$j\in\N$; this implies $| \xi_j/|\eta_j| - \eta_j/|\eta_j| | < 1/j$ which
shows that $\xi_j/|\eta_j|$ has an accumulation point $\xi_0\in\Sig_1$ with
$|\xi_0|=1$; but then $\xi_0$ is also an accumulation point of
$\eta_j/|\eta_j|$ and therefore an element of $\Sig_2$ --- a contradiction.

(ii): See \cite{D96}, proof of Th.\ 1.3.6.

Assume that (iii) does not hold; for each $k$ choose conic neighborhoods $W_j^k$
($j=1,2$)  in $\R^n\setminus 0$ with
the following property: $\forall\eta\in W_j^k$ the projection $\eta/|\eta|$ to 
$S^{n-1}$ has distance less than $1/k$ to the compact set $\Ga_j\cap S^{n-1}$.

By assumption we can choose $\eta_j^k\in W_j^k$ such that $\eta_1^k + \eta_2^k 
\not\in W$. In particular, $\Ga = \R^n \setminus W$ is a nonempty  
closed cone in $\R^n$ and therefore also $\eta^k = (\eta_1^k + \eta_2^k)/\be_k$ 
with $\be_k = |\eta_1^k| + |\eta_2^k|$ is contained in $\Ga$.

Now we choose $\xi_{j,0}^k\in\Ga_j\cap S^{n-1}$ such that 
$|\xi_{j,0}^k - \eta_j^k/|\eta_j^k|| < 1/k$ and set
\[
  \xi^k = \underbrace{\frac{|\eta_1^k|}{\be_k} \xi_{1,0}^k}_{\xi_1^k\in\Ga_1} 
           + 
          \underbrace{\frac{|\eta_2^k|}{\be_k} \xi_{2,0}^k}_{\xi_2^k\in\Ga_2}
    \in \Ga_1 + \Ga_2 \; . 
\]
By (ii) of the current lemma $(\Ga_1\cup \{0\}) + (\Ga_2\cup \{0\})$ is a closed 
cone in $\R^n$; it intersects $\Ga$ only in $0$, so we can apply (i): for
all $k\in\N$ we have for some $\al > 0$
\[
  0 < \al |\xi^k| \leq |\xi^k-\eta^k| =
  \frac{1}{\be_k}\big| |\eta_1^k|(\xi_{1,0}^k-\frac{\eta_1^k}{|\eta_1^k|})
   + |\eta_2^k|(\xi_{2,0}^k-\frac{\eta_2^k}{|\eta_2^k|})  \big| < \frac{1}{k}
   \; .
\]
Sending $k\to\infty$ we conclude that $\xi^k\to 0$. By construction
the summands of $\xi^k$ are bounded: $|\xi_j^k|\leq 1$ ($j=1,2$). If $\xi_1^k$ 
would tend to $0$ then so would $\xi_2^k = \xi^k - \xi_1^k$. But since
$\xi_{j,0}^k$ are normalized this would imply that both $|\eta_j^k|/\be_k$
($j=1,2$) tend to zero yielding the contradiction 
$ 1 = (|\eta_1^k|+|\eta_2^k|)/\be_k \to 0$. Therefore the norms of 
(suitable subsequences of) $\xi_j^k$
are bounded away from zero and above. There are subsequences $\xi_{j}^{k_l}$
($l\in\N$) such that $\xi_{j}^{k_l}\to\zeta_j\not=0$; then
$\Ga_1+\Ga_2 \ni \zeta_1+\zeta_2 = 0$ and therefore $0\in \Ga_1+\Ga_2$ ---
a contradiction.  
\end{proof}

An essential new feature of microlocal analysis in the Colombeau setting is
the precise quantification of decrease properties in terms of powers of the 
regularization parameter $\eps$. Recall from \cite{GH99}, Def.\ 17 that $R\in 
\G_\tau$ is called rapidly decreasing in a cone $\Gamma$ if $\exists N$ $\forall 
p\in \N_0$ $\exists M\in \N_0$ $\forall \phi\in {\cal A}_M$ $\exists c>0$
$\exists \eta>0$:
\begin{equation}\label{rapdec}
|R(\phi_\eps,x)| \le c\eps^{-N}(1+|x|)^{-p} \qquad (x\in \Gamma,\, 0<\eps<\eta)\,.
\end{equation}
The following lemma shows that on closed cones in the complement of the cone of
irregular directions of $U\in \G_c$, the order $N$ in (\ref{rapdec}) of rapid decrease
of the Fourier transform of $U$ can be chosen uniformly.

\begin{lemma} {\it \label{unifrapdec}
Let $U\in \G_c(\Omega)$ and let $\Gamma$ be a closed cone in the complement of $\Sigma_g(U)$.
Then $\exists N$ $\forall 
p\in \N_0$ $\exists M\in \N_0$ $\forall \phi\in {\cal A}_M$ $\exists c>0$
$\exists \eps_0>0$:
\begin{equation}\label{unifrapdecequ}
|{\cal F}(U(\phi_\eps,\,.\,))(\xi)| \le c\eps^{-N}(1+|\xi|)^{-p} \qquad 
(\xi\in \Gamma,\, 0<\eps<\eps_0)\,.
\end{equation}
}
\end{lemma}
\begin{proof}
For any $\eta\in \Ga$ there exists an open conic neighborhood $\Ga(\eta)$ such that
$\exists N(\eta)$ $\forall 
p\in \N_0$ $\exists M(\eta,p)\in \N_0$ $\forall \phi\in {\cal A}_M$ $\exists c(\eta,p,\phi)>0$
$\exists \eps_0(\eta,p,\phi)>0$ such that (\ref{unifrapdecequ}) holds with this set of constants
on $\Ga(\eta)$. The sets $\Ga(\eta)\cap S^{n-1}$ are open in $S^{n-1}$ and form a covering of
the compact set $\Ga\cap S^{n-1}$. Thus there exist $\eta_1,\dots,\eta_m \in \Ga$ such that
$$
\Ga\cap S^{n-1} \subseteq \bigcup_{j=1}^m (\Ga(\eta_j)\cap S^{n-1})\,.
$$
Consequently, $\Ga$ is contained in the union of the $\Ga(\eta_j)$ ($1\le j\le m$). Now set
$N = \max_{1\le j \le m} N(\eta_j)$ to finish the proof.  
\end{proof}

Following the terminology of \cite{MO} we will say that the wave front sets of two elements
$V_1$, $V_2$ of $\G$ are {\it in favorable position} if $\WF_g(V_1) + \WF_g(V_2)$ does not contain any zero 
direction (i.e.\ any element of the form $(x,0)$). For $V_1=\iota(v_1)$, $V_2=\iota(v_2)$ distributions 
(in which case $\WF(v_i)$ and $\WF_g(\iota(v_i))$ coincide by \cite{NPS98}, Th.\ 3.8 and 
\cite{GH99}, Cor.\ 24) this condition ensures that the Fourier 
product $v_1 v_2$ of $v_1$ and $v_2$ exists in $\D'$ (\cite{MO}, Prop.\ 6.3). Also, by \cite{MO}, Prop.\
10.3, $V_1 V_2$ is associated with $v_1 v_2$ in this case. Moreover, the wave front sets of $v_1$, $v_2$ and
$v_1 v_1$ are related by (see \cite{H}, Th.\ 8.2.10)
\begin{equation}\label{distwfrel}
\WF(v_1 v_2) \subseteq (\WF(v_1)+\WF(v_2))\cup \WF(v_1)\cup \WF(v_2) \,.  
\end{equation}
Our aim in the remainder of this section is to prove the analog of relation (\ref{distwfrel}) for elements 
of $\G$ whose generalized wave front sets are in favorable position, where the product is to be taken in 
the algebra $\G$. In the following section it will turn out that the inclusion will in general break down
if the assumption of a favorable position of the wave front sets is dropped.

\begin{prop}{\it \label{prodlemma}
Let $V_1$, $V_2$ $\in\G_c(\R^n)$ and suppose that $0\not \in \Sig_g(V_1)+\Sig_g(V_2)$. Then
\begin{equation}\label{sig_lem}
  \Sig_g (V_1 V_2) \subseteq  \ovl{\Sig_g(V_1)+\Sig_g(V_2)}^{\R^n\setminus 0} = (\Sig_g(V_1)+\Sig_g(V_2))
\cup \Sig_g(V_1) \cup \Sig_g(V_2)\, .
\end{equation}
}
\end{prop}
\begin{proof} 
Choose representatives $v_1$, $v_2$ with compact support; with the short hand
notation $w_j^\eps(\eta) = \F(v_j(\phi_\eps,.))(\eta)$ we have to estimate 
\[ 
  (2\pi)^n \F\big( v_1(\phi_\eps,.)v_2(\phi_\eps,.) \big)(\xi) = 
   w_1^\eps * w_2^\eps\,
   (\xi) = \int\limits_{\R^n} w_1^\eps(\xi-\eta) w_2^\eps(\eta)\, d\eta 
 \]
in a suitable conic neighborhood of any point $\xi_0$ in the complement of the 
right hand side of (\ref{sig_lem}). 

Let $\Ga_0$ be an open cone containing $\ovl{\Sig_g(V_1)+\Sig_g(V_2)}^{\R^n
\setminus 0}$ such that $\xi_0 \not\in \ovl{\Ga_0}$. By Lemma \ref{joshi} (iii)
there exist open cones $\Ga_j \supseteq \Sig_g(V_j)$ ($j=1,2$) such that
$\Ga_1+\Ga_2 \subseteq \Ga_0$. Further, we set $\Ga = \R^n\setminus \ovl{\Ga_0}$.
We claim that $w_1^\eps * w_2^\eps$ is rapidly decreasing in $\Ga$.

To show this we write
$$
w_1^\eps * w_2^\eps(\xi) = 
\underbrace{\int_{\Ga_2^c} w_1^\eps(\xi-\eta)
w_2^\eps(\eta)\,d\eta}_{I^\eps_1(\xi)}
+\underbrace{\int_{\Ga_2} w_1^\eps(\xi-\eta)
w_2^\eps(\eta)\,d\eta}_{I^\eps_2(\xi)}
$$
and estimate the summands individually.

Substituting $\eta' = \xi -\eta$, $I_1^\eps$ takes the form
$$
 I_1^\eps(\xi) =
\underbrace{\int\limits_{(\{\xi\} - \Ga_2^c)\cap \Ga_1} w_2^\eps(\xi-\eta') w_1^\eps(\eta')\, 
d\eta'}_{I_{11}^\eps(\xi)}
+ \underbrace{\int\limits_{(\{\xi\} - \Ga_2^c)\cap \Ga_1^c} w_2^\eps(\xi-\eta') w_1^\eps(\eta')\, 
d\eta'}_{I_{12}^\eps(\xi)}
$$
\fbox{$I_{12}^\eps$} By Lemma \ref{unifrapdec} $\exists N$ $\forall 
p\in \N_0$ $\exists M\in \N_0$ $\forall \phi\in {\cal A}_M$ $\exists c>0$
$\exists \eps_0>0$:
$$
|w_2^\eps(\xi-\eta')| \le c(1+|\xi-\eta'|^2)^{-p}\eps^{-N} \qquad (\eta' \in \{\xi\}-\Ga_2^c,\,
\eps \in (0,\eps_0))
$$ 
and $\exists N'$ $\forall 
p'\in \N_0$ $\exists M'\in \N_0$ $\forall \phi\in {\cal A}_{M'}$ $\exists c'>0$
$\exists \eps_0'>0$: 
$$
|w_1^\eps(\eta')| \le c'(1+|\eta'|^2)^{-p'}\eps^{-N'} \qquad (\eta' \in \Ga_1^c,\,
\eps \in (0,\eps_0'))
$$  
Thus by Peetre's inequality we obtain (for $\eps$ small and $\phi\in {\cal A}_{\max(M,M')}$):
$$
|I_{12}^\eps(\xi)| \le c''\eps^{-N-N'} (1+|\xi|^2)^{-p}\int_{\R^n} (1+|\eta'|^2)^{p-p'} \,d\eta'\,.   
$$
This last integral is convergent for $p'>p+\frac{n}{2}$, so $I_{12}^\eps$ is rapidly decreasing in $\Ga$.

\fbox{$I_{11}^\eps$} We abbreviate the domain of integration by $B_\xi = \Ga_1\cap 
(\{\xi\}-\Ga_2^c)$. For $\eta'\in B_\xi$, $w_1^\eps$ is tempered in $\eta'$ and $w_2^\eps$
is rapidly decreasing in $\xi-\eta'$. (For later use we note here that since $\Ga_1 \subseteq
(\Ga-\Ga_2)^c$, the same decrease properties for $w_1^\eps$ and $w_2^\eps$ in fact hold on all of 
$\Ga_1$. The following estimates thus remain valid upon replacing $B_\xi$ by $\Ga_1$.)
Hence
$$
|I_{11}^\eps(\xi)|\le c\eps^{-N} \int_{B_\xi}(1+|\xi-\eta'|)^{-p}(1+|\eta'|)^M\,d\eta'\,.
$$
Supposing $|\xi|\ge 1$ and setting $\xi_0=\frac{\xi}{|\xi|}$ this equals
\begin{eqnarray}
&& c\eps^{-N} |\xi|^{M-p} \int_{B_\xi} \left(\frac{1}{|\xi|} + 
\left|\xi_0-\frac{\eta'}{|\xi|}\right|\right)^{-p}
\left(\frac{1}{|\xi|} + \frac{\eta'}{|\xi|}\right)^M\,d\eta' \nonumber\\ 
&& \le c\eps^{-N} |\xi|^{M+n-p} \int_{\frac{1}{|\xi|} B_\xi} \left(\frac{1}{|\xi|} + 
\left|\xi_0-\eta\right|\right)^{-p}
\left(1 + |\eta|\right)^M\,d\eta \label{upint}
\end{eqnarray}
Since $\ovl{\Ga}\cap\ovl{\Ga}_1 = \{0\}$, by (i) of Lemma \ref{joshi} we have
\begin{eqnarray}
 \exists \alpha>0: &\quad |\xi_0 - \eta| \ge \alpha |\eta| \qquad 
&\forall \xi_0 \in S^{n-1}\cap \Ga, \, \forall \eta \in \Ga_1 \label{21}\\
 \exists \beta>0: &\quad |\xi_0 - \eta| \ge \beta |\xi_0| = \beta\qquad 
&\forall \xi_0 \in S^{n-1}\cap \Ga, \, \forall \eta \in \Ga_1\label{22}
\end{eqnarray}
We now split the domain of integration in (\ref{upint}) into the parts $B_1 =
\frac{1}{|\xi|} B_\xi \cap \{|\eta|\le \frac{1}{2}\}$ and $B_2 =
\frac{1}{|\xi|} B_\xi \cap \{|\eta|> \frac{1}{2}\}$. Then by (\ref{22})
$$
\int_{B_1} \left(\frac{1}{|\xi|} + \left|\xi_0-\eta\right|\right)^{-p}
\left(1 + |\eta|\right)^M\,d\eta \le  \beta^{-p} \int_{B_1} (1+|\eta|)^M\,d\eta
\le const\,.
$$
Also, by (\ref{21}),
$$
\int_{B_2} \left(\frac{1}{|\xi|} + \left|\xi_0-\eta\right|\right)^{-p}
\left(1 + |\eta|\right)^M\,d\eta \le  
\alpha^{-p} \int_{B_2} |\eta|^{-p} (1+|\eta|)^M\,d\eta \le const
$$
for $p>M-n$. It follows that $I_{11}^\eps$ and hence also 
$I_{1}^\eps$ is rapidly decreasing in $\Ga$.

Turning now to $I_2^\eps$, we first note that $w_2^\eps$ is tempered on the domain
of integration. Moreover, since $\Ga_2\subseteq (\Ga - \Ga_1)^c$ it follows that
$w_1^\eps$ is rapidly decreasing in $\xi-\eta$ in said domain (again by Lemma 
\ref{unifrapdec}). Thus the same reasoning as in the case of $I_{11}^\eps$ 
(cf.\ the above remark) shows
that $I_2^\eps$ is rapidly decreasing in $\Ga$ as well, which completes the proof.
\end{proof}

\begin{thm}{\it \label{wfth}
Let $U_1$, $U_2$ be elements of $\G(\Omega)$ whose wave front sets are in favorable position.
Then
\begin{equation}\label{wf_thm}
    \WF_g (U_1 U_2) \subseteq  (\WF_g(U_1)+\WF_g(U_2))  \cup 
    \WF_g(U_1) \cup \WF_g(U_2)\;.
\end{equation}
}
\end{thm}
\begin{proof}

Let $(x,\xi) \not \in r.h.s.$. Then for any $\vphi \in \D(\Omega)$ with 
$\vphi(x)\not=0$ and support sufficiently close to $x$ we have
$\xi \not \in \Sig_g(\vphi U_i)$ ($i=1,2$). Also, by Lemma 
\ref{joshi} (ii) and (iii), since $\xi$ is not contained in
$$
\ovl{\Sig_{g,x}(U_1) + \Sig_{g,x}(U_2)}^{\R^n\setminus 0}
$$
there exist open conic neighborhoods $\Ga_i$ of $\Sig_{g,x}(U_i)$ in 
$\R^n\setminus 0$ such that
$\xi \not\in \Ga_1 + \Ga_2$ and $0\not\in \Ga_1 +\Ga_2$. 
Thus, by \cite{GH99}, (13), if  the support of $\vphi$ is close enough to $x$ 
we also have 
$$
\xi \not \in \Sig_g(\vphi U_1) + \Sig_g(\vphi U_2) \subseteq \Ga_1 + \Ga_2\,.
$$
Since $0\not\in \Sig_g(\vphi U_1) + \Sig_g(\vphi U_2)$, 
$$
\xi \not \in \ovl{\Sig_{g}(\vphi U_1) + \Sig_{g}(\vphi U_2)}^{\R^n\setminus 0}
= (\Sig_{g}(\vphi U_1) + \Sig_{g}(\vphi U_2))\cup \Sig_{g}(\vphi U_1) 
\cup \Sig_{g}(\vphi U_2) \,.
$$
Thus by Proposition 
\ref{prodlemma}, $\xi\not\in \Sig_g(\vphi^2 U_1 U_2)$. Again from \cite{GH99},
(13) the claim follows. \end{proof}

\section{Examples} \label{exsec}
In the previous section we have extended the validity of the wave front inclusion relation
(\ref{distwfrel}) to the product in the algebra $\G$, provided that the generalized wave front
sets of the factors are in favorable position. Contrary to the distributional situation, however,
a favorable position of the wave front sets is of course not a prerequisite for forming the 
product in the algebra. Thus the question arises whether a further extension of the classical
result to arbitrary products in $\G$ is possible. The second example in this section will 
demonstrate that this is not the case. Before we turn to this matter, we first give an
example illustrating some genuinely non-distributional effects in the application of
Theorem \ref{wfth}.
\begin{ex} \label{wfex1}
Denote by $U$ be the class in $\G(\R^2)$ of $U(\phi^{(2)},x,y)$ $=$ $\frac{1}{d(\phi)}\phi
\left(\frac{x}{d(\phi)}\right.$ $\left.- \frac{y}{\sqrt{d(\phi)}}\right)$. 
As the results of the following 
calculations are independent of the concrete value of $d(\phi)$ we will for simplicity 
assume that $d(\phi)=1$ and we will abbreviate $U(\phi^{(2)}_\eps,x,y)$ by $u_\eps(x,y)
= \frac{1}{\eps}\phi \left(\frac{x}{\eps} - \frac{y}{\sqrt{\eps}}\right)$. 
As a matter of fact, this assumption effectively transfers the problem into the setting
of the special Colombeau algebra. It is easily seen that $U \approx \delta(x)\otimes
1(y)$.

Further, let $A = \iota(\frac{1}{x+i0})
\in \G(\R)$ and define $B\in \G(\R^2)$ by $B(\phi^{(2)},x,y) = A(\phi,\sqrt{d(\phi)}x+y)$.
Employing the same simplification as above we will write $b_\eps(x,y)$ for
$$
B(\phi^{(2)}_\eps,x,y) = \int_0^\infty \frac{\phi_\eps(\sqrt{\eps}x+y-z)
-\phi_\eps(\sqrt{\eps}x + y +z)}{z}\,dz - i\pi\phi_\eps(\sqrt{\eps}x + y)
$$
Let us first determine $\WF_g(U)$. To begin with, we claim that $\supp(U) = \{0\}\times \R$.
Indeed, the inclusion $\subseteq$ is obvious. Conversely, let $(0,a)\in \{0\}\times \R$ and
set $x_\eps = \eps+a\sqrt{\eps}$, $y_\eps = \sqrt{\eps} + a$. Then $(x_\eps,y_\eps)$ is the
representative of a compactly supported generalized point (cf.\ \cite{OK99}) supported in 
any ball $B_r((0,a))$ ($r>0$ arbitrary). Since $u_\eps(x_\eps,y_\eps) = \frac{1}{\eps}\phi(0)$,
the claim follows from \cite{OK99}, Th.\ 2.4. To determine an upper bound for $\WF_g(U)$ we
note that setting $P_\eps = \partial_y + \sqrt{\eps} \partial_x$ we have $P_\eps u_\eps = 0$.
Thus by (\ref{prop}) the set of characteristic directions of $P$ provides such an upper bound.
By Example \ref{firstex}, $\mbox{Char}(P)\subseteq \R^2\times \{(\xi,0)\mid \xi\not=0\}$.
Next, we show that $\{(\xi,0)\mid \xi\not=0\}\subseteq \Sigma_{g,(0,a)}(U)$ for any
$a\in \R$. To this end, let $\psi(x,y) = f(x)g(y) \in \D(\R^2)$, $f(x) \equiv 1$, 
$g(y)\equiv 1$ near $x=0$ resp.\ $y=0$, and $g$ positive. For $\eps$ sufficiently small, 
$\{x\mid g(y)\phi_\eps(x-\sqrt{\eps}y)\not=0\} \subseteq \{x\mid f(x)=1\}$, so
$\psi u = g u$. Hence
\begin{eqnarray*}
&&{\cal F}(\psi u_\eps)(\xi,0) = \int e^{-ix\xi}\phi_\eps(x-\sqrt{\eps}y)g(y)\,dxdy 
=\hat{\phi}_\eps(\xi)\hat{g}(\sqrt{\eps}\xi) \\
&& = \hat{\phi}(\eps\xi)\hat{g}(\sqrt{\eps}\xi) \to \hat{\phi}(0)\hat{g}(0)
= \hat{g}(0) \not=0 \qquad (\eps\to 0)\,. 
\end{eqnarray*}
which shows that ${\cal F}(\psi u_\eps)$ is not rapidly decreasing in the direction
$(\xi,0)$. Replacing $g$ by $\tau_a g = g(\,.\,-a)$ and setting $\psi_a(x,y) =
f(x)\tau_a g(y)$ we obtain ${\cal F}(\psi_a u_\eps) = \hat{\phi}(\eps\xi)
\hat{g}(\sqrt{\eps}\xi) e^{-ia\xi}$, so the same reasoning gives 
$\{(\xi,0)\mid \xi\not=0\}$ $\subseteq$ $\Sigma_{g,(0,a)}(U)$. Summing up, we have
shown
\begin{equation}\label{wfu}
\WF_g(U) = \{0\} \times \R \times \R\setminus 0 \times \{0\} \,.  
\end{equation}
Turning now to $B$, we first show that
$\singsupp(B)$ $=$ $\{(0,0)\}$. Let $(x,y) \in K\subset\subset \R^2\setminus 0$.
Since $\phi(\frac{x}{\sqrt{\eps}} + \frac{y\pm z}{\eps})
\equiv 0$ near $z=0$ for $\eps$ small, we can write
\begin{eqnarray*}
b_\eps(x,y) &=& \frac{1}{\eps}\int_0^\infty 
\frac{\phi(\frac{x}{\sqrt{\eps}} + \frac{y -z}{\eps})}{z}\,dz  
- \frac{1}{\eps}\int_0^\infty 
\frac{\phi(\frac{x}{\sqrt{\eps}} + \frac{y+z}{\eps})}{z}\,dz\\
&=& \left\{ 
\begin{array}{ll}
\hphantom{-} \frac{1}{\eps} \int_{-\infty}^\infty 
\frac{\phi(s)\,ds}{y+\sqrt{\eps}x-\eps s} 
& \quad \frac{x}{\sqrt{\eps}} + \frac{y}{\eps} > d(\phi) \\
-\frac{1}{\eps} \int_{-\infty}^\infty \frac{\phi(s)\,ds}{\eps s-y -\sqrt{\eps}x} 
& \quad \frac{x}{\sqrt{\eps}} + \frac{y}{\eps} < -d(\phi) 
\end{array}
\right.
\end{eqnarray*}
In any case, 
$$
\partial_y^k \partial_x^l(b_\eps) = \pm\frac{1}{\eps}\eps^{\frac{l}{2}}
\int_{-\infty}^\infty \frac{\phi(t)\,dt}{(-\eps t \pm (y+\sqrt{\eps}x))^{k+l+1}}
= {\cal O}(\frac{1}{\eps})\,,
$$ so $B$ is an element of $\G^\infty$ off $(x,y)=(0,0)$ (but clearly not in any
neighborhood of $(0,0)$ itself).
Setting $P_\eps = \partial_x - \sqrt{\eps} \partial_y$ we have $P_\eps b_\eps = 0$,
so by the same reasoning as above we have $\WF_g(B)\subseteq \mbox{Char}_g(P) =
\R^2 \times \{\xi=0\}$. To determine $\Sigma_{g,(0,0)}B$ it therefore remains to 
estimate (with $\psi$ as above)
\begin{eqnarray*}
&&{\cal F}(\psi b_\eps)(0,\eta) = \int e^{-iy\eta} \int_0^\infty 
\frac{\phi_\eps(\sqrt{\eps}x+y-z) 
-\phi_\eps(\sqrt{\eps}x + y +z)}{z}\,dz \cdot \\ 
&& \cdot f(x)g(y)\,dx dy - i\pi\int e^{-iy\eta} \phi_\eps(\sqrt{\eps}x + y) 
f(x)g(y)\,dx dy \\
&& = \frac{1}{2\pi} \int f(x)
\underbrace{({\cal F}_{y\to \eta'}(((\mbox{vp}(\frac{1}{x})-i\pi\delta)
*\phi_\eps)(\sqrt{\eps}x+y)*\hat g)(\eta)}_{=: I^\eps(x,\eta)}\,dx
\end{eqnarray*}
Using ${\cal F}(\mbox{vp}(\frac{1}{x})-i\pi\delta) = c H$ (with $H$ the Heaviside function) 
we have
$$
I^\eps(x,\eta) = c\int_0^\infty \hat g(\eta-\eta') e^{i\eta'\sqrt{\eps}x} 
\hat\phi(\eps\eta')\, d\eta'
$$
Suppose now that $\eta<0$. Then in the domain of integration of $I^\eps$,
$|\eta-\eta'| \ge |\eta|$. Thus for any $l\in \N$ we get
$$
|\hat g(\eta-\eta')| \le c_l(1+|\eta-\eta'|)^{-2l} \le
c_l(1+|\eta|)^{-l}(1+|\eta-\eta'|)^{-l}\,.
$$
But then for $l$ sufficiently large $|I^\eps(x,\eta)| \le c'(1+|\eta|)^{-l}$, 
implying $\Sigma_{g,(0,0)}(B) \subseteq \{(0,\eta)\mid \eta>0\}$. Since we have seen above that
$\Sig_{g,(0,0)} \not= \emptyset$, this implies
\begin{equation}\label{wfb}
\WF_g(B) = \{(0,0)\} \times \{0\}\times \R^+ \,.  
\end{equation}
By (\ref{wfu}) and (\ref{wfb}), the wave front sets of $U$ and $B$ are in 
favorable position, so Theorem \ref{wfth} gives 
\begin{equation}\label{wfbuincl}
  \Sigma_{g,(0,0)}(BU)\subseteq \R\times \R^+\,.
\end{equation}
We are now going to establish also the inverse inclusion to (\ref{wfbuincl}). With
$\psi$ as above, we have to analyze 
$$
{\cal F}(\psi b_\eps u_\eps)(\xi,\eta) =
{\cal F}_{x\to\xi}(\underbrace{{\cal F}_{y\to\eta}(g(y)b_\eps(x,y)u_\eps(x,y))
}_{=: J^\eps(x,\eta)}f(x))\,.
$$
Here, $(2\pi)^2 J^\eps(x,\eta) = (\hat g * {\cal F}_{y\to\eta'}(b_\eps(x,y))
* {\cal F}_{y\to\eta'}(u_\eps(x,y)))(\eta)$, and a short calculation gives:
\begin{eqnarray*}
&&{\cal F}_{y\to\eta'}(b_\eps(x,y))(\eta') = -2\pi i e^{i\eta'\sqrt{\eps}x}H(\eta')
\hat\phi(\eps\eta')\\  
&& {\cal F}_{y\to\eta'}(u_\eps(x,y))(\eta') = \frac{1}{\sqrt{\eps}}  e^{-i\eta'\frac{x}{\sqrt{\eps}}}
\hat\phi(-\sqrt{\eps}\eta')
\end{eqnarray*}
Inserting this and substituting $x' = \frac{x}{\sqrt{\eps}}$ we obtain
\begin{eqnarray}\label{ftbu}
&{\cal F}(\psi b_\eps u_\eps)(\xi,\eta)=
\frac{1}{2\pi i} \int\int\int_0^\infty e^{-i\xi\sqrt{\eps}x' +ix'(\eps\eta''-\eta')}
\hat g(\eta-\eta'-\eta'')\cdot& \nonumber\\ 
& \cdot\hat\phi(\eps\eta'')\hat\phi(-\sqrt{\eps}\eta')
\,d\eta'' d\eta' f(\sqrt{\eps} x')\, dx'&  
\end{eqnarray}
For $\eps\to 0$ this converges to
\begin{eqnarray*}
&&\hphantom{=}\frac{1}{2\pi i} \int\int\int_0^\infty e^{-i x'\eta'}\hat g(\eta-\eta'-\eta'')
\, d\eta'' d\eta' dx'\\
&&= \frac{1}{2\pi i} \int \int_{-\infty}^z\hat g(\xi)\,d\xi\int e^{-ix'(\eta-z)}\,
dx' dz = \frac{1}{2\pi i} \int_{-\infty}^\eta \hat g(\xi)\, d\xi\,.
\end{eqnarray*}
It follows that for $\eta>0$ and any $\xi$, ${\cal F}(\psi b_\eps u_\eps)$ is not
rapidly deceasing in the direction $(\xi,\eta)$. Thus, in fact
\begin{equation}\label{buwf}
  \Sigma_{g,(0,0)}(BU) = \R\times \R^+\,.
\end{equation}
\end{ex}

\begin{ex} \label{moex}
Let $U$ as in Example \ref{wfex1} and set $v = u(\phi^{(2)},x,y) = \frac{1}{d(\phi)}\phi
\left(\frac{x}{d(\phi)}\right.$ $-$ $\left. \frac{y}{\sqrt{d(\phi)}}\right)$. 
Employing the same notational
simplifications as in the previous example we have $v_\eps(x,y)$
$=$ $\frac{1}{\eps}\phi \left(\frac{x}{\eps} + \frac{y}{\sqrt{\eps}}\right)$.
Again, $\supp(V) = \{0\}\times \R$, $V \approx \delta(x)\otimes 1(y)$, and
\begin{equation}\label{wfv}
\WF_g(V) = \{0\} \times \R \times \R\setminus 0 \times \{0\} \,.  
\end{equation}
Thus the wave front sets of $U$ and $V$ are not in favorable position and we shall
demonstrate that in fact the conclusion of Theorem \ref{wfth} is violated for the product
$W = UV$.

We first show that $\supp(W) = \{(0,0)\}$. To this end we again utilize \cite{OK99}, 
Th.\ 2.4. We only have to show that $W$ vanishes in a suitable neighborhood of any
point $(0,a)$ with $a\not=0$. We shall assume $a>0$ (the other case being analogous)
and we choose some $r>0$ such that $B_r((0,a))$ does not contain $(0,0)$. Now let
$(x_\eps,y_\eps)$ be a representative of any generalized point supported in 
$B_r((0,a))$. Then if $x_\eps\ge 0$  and $\eps$ is sufficiently small we have 
$$
\frac{x_\eps}{\eps} + \frac{y_\eps}{\sqrt{\eps}} \ge \frac{a-r}{\sqrt{\eps}}
$$
and similarly, for $x_\eps \le 0$, 
$\frac{x_\eps}{\eps} - \frac{y_\eps}{\sqrt{\eps}} \le -\frac{a-r}{\sqrt{\eps}}$.
Thus, for small $\eps$, one of the factors of $w_\eps$ always vanishes. This means
that $W$ vanishes on all compactly supported points in $B_r((0,a))$, so $W=0$ on
$B_r((0,a))$.

It remains to determine $\Sigma_{g,(0,0)}(W)$, to which end we choose $f$, $g$ as
in Example \ref{wfex1}. As above, $f^2(x)g^2(y)u_\eps(x,y)v_\eps(x,y) =
g^2(y)u_\eps(x,y)v_\eps(x,y)$ for $\eps$ small. Thus it suffices to consider
\begin{eqnarray*}
&& {\cal F}(gu_\eps gv_\eps)(\xi,\eta) = 
\int\int \hat\phi_\eps(\xi')\hat\phi_\eps(\xi-\xi') \hat g(\eta'+\sqrt{\eps}\xi')
\hat g(\eta-\eta' \\
&&-\sqrt{\eps}(\xi-\xi'))\,d\xi' d\eta' 
 = \int \hat\phi(\eps\xi')\hat\phi(\eps(\xi-\xi')) \hat g * \hat g
(\eta+2\sqrt{\eps}\xi'-\sqrt{\eps}\xi)\, d\xi'  \\
&& = \frac{1}{\sqrt{\eps}}\int \hat\phi(\sqrt{\eps}\xi'')\hat\phi(\eps\xi-\sqrt{\eps}\xi'') 
\hat g * \hat g(\eta+2\xi''-\sqrt{\eps}\xi)\, d\xi''  
\end{eqnarray*}
The integral in this equation converges to $\int \hat g * \hat g(2\xi'')\,d\xi''$, 
so ${\cal F}(gu_\eps gv_\eps)(\xi,\eta) \sim {\cal O}(\frac{1}{\sqrt{\eps}})$. In 
particular, ${\cal F}(gu_\eps gv_\eps)$ is not rapidly decreasing in any direction
$(\xi,\eta)$. Thus
\begin{equation}\label{moexequ}
\begin{array}{rcl}
\WF_g(UV) &=& \{(0,0)\}\times \R^2\setminus 0 \not\subseteq 
\{0\} \times \R \times \R\setminus 0 \times \{0\} \\ 
&=& (\WF_g(U) + \WF_g(V)) \cup \WF_g(U)\cup \WF_g(V)   \,.
  \end{array}
\end{equation}
\end{ex}

\section{Consequences for 
microlocal properties of pullbacks}

In this final section we are going to compare our previous considerations with
an alternative (classical) approach to products of distributions and their microlocal 
properties (\cite{H}, section 8.2). In this approach, one
considers the product of two distributions $u$ and $v$ (provided it exists) as the restriction of 
their tensor product to the diagonal, i.e., as the
pullback of $u\otimes v$ under $d: x \to (x,x)$. The basic properties of this 
operation then follow directly from the general theorem about microlocal transformation 
under composition with smooth maps (\cite{H}, Th.\ 8.2.4)

As a first step, we note that the tensor product of Colombeau generalized functions
is ``well behaved'' from a microlocal point of view.  Recall that for open subsets $\Omega
\subseteq\R^m$, $\Omega'\subseteq\R^n$ and $U\in\G(\Omega)$, $V\in\G(\Omega')$, 
$U\otimes V$ is represented by $(\phi^{(m)}\otimes\phi^{(n)},x,y) \to
u(\phi^{(m)},x) v(\phi^{(n)},y)$.

\begin{lemma} {\it \label{tensorwf}
\begin{eqnarray}
&&  \WF_g(U\otimes V) \subseteq \Big(\WF_g(U)\boxtimes\WF_g(V)\Big)
    \cup \Big(\big(\supp U\times \{0\}\big) \boxtimes \WF_g(V)\Big) \nonumber\\
&&   \hphantom{\WF_g(U\otimes V) \subseteq}     
     \cup\Big( \WF_g(U)\boxtimes \big(\supp V \times \{0\}\big)\Big)
\end{eqnarray}
where $\Gamma_1 \boxtimes \Gamma_2 := \{ (x,y,\xi,\eta) \mid
(x,\xi)\in\Gamma_1, (y,\eta)\in\Gamma_2\}$ for arbitrary subsets
$\Gamma_1\subseteq \Omega\times\R^m$, $\Gamma_2\subseteq \Omega'\times\R^n$.
}
\end{lemma}
\begin{proof}
This is a straightforward adaptation of the proof of the corresponding distributional result
(see e.g. \cite{FJ}, Th.\ 11.2.1).
\end{proof}

If $f: \Omega_1 \map \Omega_2$ is a smooth map between open subsets $\Omega_i \subseteq \R^{n_i}$ 
and $u\in \D'(\Omega_2)$ then the classical condition ensuring existence of $f^* u$ is
\[
  \WF(u) \cap N_f = \emptyset
\]
where $N_f = \{(f(x_1),\xi_2)\in \Omega_2\times\R^{n_2} 
\mid {}^t{f'(x_1)}\xi_2 = 0\}$ (cf.~\cite{H}, Th.\ 8.2.4). Furthermore, in
this case
\begin{equation} \label{pullbackwf}
  \WF(f^* u) \subseteq f^* \WF(u)
\end{equation}
where $f^*\Gamma = \{ (x_1, {}^t f'(x_1) \xi_2 \mid
(f(x_1),\xi_2)\in\Gamma \}$ for $\Gamma\subseteq \Omega_2\times\R^{n_2}$. 

Just as the product of generalized functions can be carried out unrestrictedly in the 
Colombeau algebra, it is also possible to form pullbacks of Colombeau generalized functions
under arbitrary smooth maps: For $U\in \G$ and $f$ smooth, $U\circ f$ is defined by
componentwise composition. Moreover, for $U$, $V$ in $\G(\Omega)$ we have 
$U V = (U\otimes V)\circ d$ with $d: \Omega\map \Omega\times \Omega$, $x \mapsto (x,x)$ (simply 
observe that $U V$ and $(U\otimes V)\circ d$ have identical representatives).
We are thus in a position to review our previous examples in this picture:

\begin{ex} \label{12}
Let $U$ and $V$ as in Example \ref{moex} and set $T = U\otimes V$. 
By Lemma \ref{tensorwf},
\[
  \WF_g(T) \subseteq \{((0,x_2,0,x_4),(\xi_1,0,\xi_3,0)) \mid x_2, x_4
  \in\R, (\xi_1,\xi_3) \in \R^2\setminus \{(0,0)\} \} \, .
\]
In fact, we have equality in the above relation since by using cutoff
functions of tensor product form the reasoning leading
to (\ref{wfu}) can be carried out in parallel in the independent factors
corresponding to $U$ and $V$. In the notation of Example \ref{moex} we have
$W = T \circ d$ and $\WF_g(W) = \{ (0,0)\} \times \R^2\setminus 0$. A simple
computation shows that $d^* \WF_g(T) = \{((0,r),(\mu,0)) \mid r,\mu\in\R \}$
which implies that
\begin{equation}\label{pullwf}
          \WF_g(d^*T)  \not\subseteq   d^* \WF_g(T) \, .
\end{equation}
\end{ex}

Note that $N_d \cap \WF_g(T) \not= \emptyset$ in Example \ref{12}. We conclude
that the validity (\ref{pullbackwf}) cannot be extended to arbitrary pullbacks 
of Colombeau functions under smooth maps. 

\begin{rem} 
We note that a common feature of the examples introduced in Section \ref{exsec} 
is that they are formed as pullbacks
of (canonical images of) distributions under {\it generalized} maps. Thus the fact that --- 
contrary to the distributional setting --- composition of generalized functions can be 
carried out in $\G$ (subject to certain growth conditions) can be viewed as one of the causes
of the new microlocal effects presented there.
\end{rem}

{\bf Acknowledgements:} We would like to thank Michael Oberguggenberger for 
suggesting the concrete form of Example \ref{moex} to us, as well as for 
several discussions that importantly contributed to the final form of the
manuscript.

\end{document}